\vsize=9.0in\voffset=1cm
\looseness=2


\message{fonts,}

\font\tenrm=cmr10
\font\ninerm=cmr9
\font\eightrm=cmr8
\font\teni=cmmi10
\font\ninei=cmmi9
\font\eighti=cmmi8
\font\ninesy=cmsy9
\font\tensy=cmsy10
\font\eightsy=cmsy8
\font\tenbf=cmbx10
\font\ninebf=cmbx9
\font\tentt=cmtt10
\font\ninett=cmtt9

\font\ninesl=cmsl9
\font\eightsl=cmsl8

\font\nineit=cmti9
\font\eightit=cmti8

\skewchar\ninei='177 \skewchar\eighti='177
\skewchar\ninesy='60 \skewchar\eightsy='60

\def\eightpoint{\def\rm{\fam0\eightrm} 
\normalbaselineskip=9pt
\normallineskiplimit=-1pt
\normallineskip=0pt

\textfont0=\eightrm \scriptfont0=\sevenrm \scriptscriptfont0=\fiverm
\textfont1=\ninei \scriptfont1=\seveni \scriptscriptfont1=\fivei
\textfont2=\ninesy \scriptfont2=\sevensy \scriptscriptfont2=\fivesy
\textfont3=\tenex \scriptfont3=\tenex \scriptscriptfont3=\tenex
\textfont\itfam=\eightit  \def\it{\fam\itfam\eightit} 
\textfont\slfam=\eightsl \def\sl{\fam\slfam\eightsl} 

\setbox\strutbox=\hbox{\vrule height6pt depth2pt width0pt}%
\normalbaselines \rm}

\def\ninepoint{\def\rm{\fam0\ninerm} 
\textfont0=\ninerm \scriptfont0=\sevenrm \scriptscriptfont0=\fiverm
\textfont1=\ninei \scriptfont1=\seveni \scriptscriptfont1=\fivei
\textfont2=\ninesy \scriptfont2=\sevensy \scriptscriptfont2=\fivesy
\textfont3=\tenex \scriptfont3=\tenex \scriptscriptfont3=\tenex
\textfont\itfam=\nineit  \def\it{\fam\itfam\nineit} 
\textfont\slfam=\ninesl \def\sl{\fam\slfam\ninesl} 
\textfont\bffam=\ninebf \scriptfont\bffam=\sevenbf
\scriptscriptfont\bffam=\fivebf \def\bf{\fam\bffam\ninebf} 
\textfont\ttfam=\ninett \def\tt{\fam\ttfam\ninett} 

\normalbaselineskip=11pt
\setbox\strutbox=\hbox{\vrule height8pt depth3pt width0pt}%
\let \smc=\sevenrm \let\big=\ninebig \normalbaselines
\parindent=1em
\rm}

\def\tenpoint{\def\rm{\fam0\tenrm} 
\textfont0=\tenrm \scriptfont0=\ninerm \scriptscriptfont0=\fiverm
\textfont1=\teni \scriptfont1=\seveni \scriptscriptfont1=\fivei
\textfont2=\tensy \scriptfont2=\sevensy \scriptscriptfont2=\fivesy
\textfont3=\tenex \scriptfont3=\tenex \scriptscriptfont3=\tenex
\textfont\itfam=\nineit  \def\it{\fam\itfam\nineit} 
\textfont\slfam=\ninesl \def\sl{\fam\slfam\ninesl} 
\textfont\bffam=\ninebf \scriptfont\bffam=\sevenbf
\scriptscriptfont\bffam=\fivebf \def\bf{\fam\bffam\tenbf} 
\textfont\ttfam=\tentt \def\tt{\fam\ttfam\tentt} 

\normalbaselineskip=11pt
\setbox\strutbox=\hbox{\vrule height8pt depth3pt width0pt}%
\let \smc=\sevenrm \let\big=\ninebig \normalbaselines
\parindent=1em
\rm}

\message{fin format jgr}

\magnification=1200
\font\Bbb=msbm10
\def\C{\hbox{\Bbb C}}
\def\R{\hbox{\Bbb R}}
\def\v{\varphi}
\def\pa{\partial}
\def\b{\backslash}
\def\ep{\varepsilon}
\setbox3=\vtop{
\hbox{Universit\'e Pierre et Marie Curie, Math\'ematiques, case 247}
\hbox{4 place Jussieu, 75252 Paris, France}
\hbox{e-mail: henkin@math.jussieu.fr}
}

\centerline{\bf Cauchy-Pompeiu type formulas for $\bar\pa$ on}
\centerline{\bf affine algebraic Riemann surfaces and some applications}

\vskip 2 mm
\centerline{\bf G.M.Henkin}
\smallskip
\centerline{\it Dedicated to Oleg Viro on the occasion of his 60th birthday}
\bigskip

{\bf Abstract}
\smallskip

We present explicit solution formulas $f=\hat R\v$ and $u=R_{\lambda}f$
for the equations $\bar\pa f=\v$ and
$(\pa+\lambda dz_1)u=f-{\cal H}_{\lambda}f$ on an
affine algebraic curve $V\subset\C^2$. Here
${\cal H}_{\lambda}f$ denotes the projection of
$f\in {\tilde W}^{1,\tilde p}_{1,0}(V)$ to the
subspace of pseudoholomorphic (1,0)-forms on $V$:
$\bar\pa {\cal H}_{\lambda}f=\bar\lambda d\bar z_1\wedge {\cal H}_{\lambda}f$.
These formulas can be interpreted as explicit versions and precisions
of the Hodge--Riemann decomposition on Riemann surfaces. The main
application consists in the construction of  the Faddeev--Green function
for $\bar\pa(\pa+\lambda dz_1)$ on $V$ as the kernel of the operator
$R_{\lambda}\circ\hat R$. This Faddeev--Green function is the main tool for
the solution of the inverse conductivity problem on bordered Riemann
surfaces $X\subset V$, that is,  for the  reconstruction
of the conductivity function
 $\sigma$ in the equation $d(\sigma d^cU)=0$ from the Dirichlet-to-Neumann
mapping $U\big|_{bX}\mapsto \sigma d^cU\big|_{bX}$.
The case $V=\C$ was treated by R.Novikov [N1]. In $\S$ 4 we give a correction
to the paper [HM], in which the case of a general algebraic curve $V$ was first
considered.
\smallskip

\smallskip

\bigskip

\noindent {\bf Introduction}
\smallskip

This paper is motivated by a problem from two-dimensional Electrical
Impedance Tomography, namely the question of how to reconstruct the
conductivity function $\sigma$ on a bordered Riemann surface
$X$ from the knowledge of the Dirichlet-to-Neumann mapping
$u\big|_{bX}\to\sigma d^cU\big|_{bX}$ for solutions $U$ of the Dirichlet
problem:
$$d(\sigma d^cU)\big|_X=0,\ \ U\big|_{bX}=u,\ \ {\rm where}\ \
d=\pa +\bar\pa,\ \ d^c=i(\bar\pa-\pa).$$
For the case $X=\Omega\subset\R^2\simeq\C$ $(z=x_1+ix_2)$ the exact
reconstruction scheme was given firstly by R.Novikov [N1] under some
restriction on the conductivity function $\sigma$.
This restriction was eliminated later by A.Nachman [Na].
\smallskip

The scheme consists in the following.
\smallskip

Let $\sigma(x)>0$ for $x\in\bar\Omega$ and $\sigma\in C^{(2)}(\bar\Omega)$.
Put $\sigma(x)=1$ for $x\in\R^2\b\bar\Omega$. The substitution
$\psi=\sqrt{\sigma}U$ transforms the equation  $d(\sigma d^cU)=0$ into the equation
$d d^c\psi={d\,d^c\sqrt{\sigma}\over \sqrt{\sigma}}\psi$ on $\R^2$.
From L.Faddeev's [F1] result (with additional arguments [BC2] and [Na])
it follows that, for each $\lambda\in\C$, there exists a unique solution
$\psi(z,\lambda)$ of the above equation,  with asymptotics
$$\psi(z,\lambda)\cdot e^{-\lambda z}
\buildrel \rm def \over =
\mu(z,\lambda)=1+o(1),\ \ z\to\infty.$$
Such a solution can be found from the integral equation
$$\mu(z,\lambda)=1+{i\over 2}\int\limits_{\xi\in\Omega}g(z-\xi,\lambda)
{\mu(\xi,\lambda)
 d\,d^c\sqrt{\sigma}\over \sqrt{\sigma}},$$
where the function
$$g(z,\lambda)={-1\over 2i(2\pi)^2}\int\limits_{w\in\C}
{e^{i(w\bar z+\bar w z)}dw\wedge d\bar w\over w(\bar w-i\lambda)},\
z\in\C,\ \lambda\in\C,$$
is called the Faddeev--Green function for the operator
$\mu\mapsto\bar\pa(\pa+\lambda dz)\mu$.
\smallskip

From work of R.Novikov [N1] it
follows that the function $\psi\big|_{b\Omega}$ can be found through
the Dirichlet-to-Neumann mapping by the integral equation
$$\eqalign{
&\psi(z,\lambda)\big|_{b\Omega}=e^{\lambda z}+\int\limits_{\xi\in b\Omega}
e^{\lambda(z-\xi)}g(z-\xi,\lambda)(\hat\Phi\psi(\xi,\lambda)-
\hat\Phi_0\psi(\xi,\lambda)),\cr
&{\rm where}\ \ \hat\Phi\psi=\bar\pa\psi\big|_{b\Omega},\ \
\hat\Phi_0\psi=\bar\pa\psi_0\big|_{b\Omega}
,\ \ \psi_0\big|_{b\Omega}=\psi,\ \
\bar\pa\pa\psi_0\big|_{\Omega}=0.\cr}$$
By results of R.Beals, R.Coifmann [BC1],
P.Grinevich, S.Novikov [GN] and R.Novikov [N2]
it then follows that $\psi(z,\lambda)$
satisfies a $\bar\pa$-equation of Bers--Vekua type with respect to
$\lambda\in\C$:
$${\pa\psi\over \pa\bar\lambda}=b(\lambda)\bar\psi,$$
where $\lambda\mapsto b(\lambda)\in L^{2+\ep}(\C)\cap L^{2-\ep}(\C)$, and
$\psi(z,\lambda)e^{-\lambda z}\to 1$ as $\lambda\to\infty$, for all $z\in\C$.
\smallskip

This $\bar\pa$-equation combined with R.Novikov's integral equation permits
us to find, starting  from the Dirichlet-to-Neumann mapping, firstly the boundary
values $\psi\big|_{b\Omega}$, secondly the "$\bar\pa$-scattering data" $b(\lambda)$,
and thirdly $\psi\big|_{\Omega}$.
\smallskip

Summarizing, the conductivity function $\sigma\big|_{\Omega}$ is thus retrieved from
the given Dirichlet-to-Neumann data by means of the scheme:
$${\rm DN\ data}\ \to\ \psi\big|_{b\Omega}
\to\  \bar\pa{\rm -scattering\ data}\ \
\to\ \psi\big|_{\Omega}\to\
{d\,d^c\sqrt{\sigma}\over \sqrt{\sigma}}\big|_{\Omega}.$$
\bigskip

{\it Main result}
\smallskip

We suppose that instead of $\C$ we have a smooth algebraic Riemann surface $V$
in $\C^2$, given by an equation $V=\{z\in\C^2;\ P(z)=0\}$, where $P$ is a holomorphic
 polynomial of degree $d\ge 1$. Put $z_1=w_1/w_0$,  $z_2=w_2/w_0$ and
suppose that the projective compactification $\tilde V$ of $V$ in
${\C}P^2\supset\C^2$ with coordinates $w=(w_0:w_1:w_2)$ intersects
$CP^1_{\infty}=\{z\in {\C}P^2;\ w_0=0\}$ transversally in $d$ points. In
order to extend the Novikov reconstruction scheme on the Riemann surface
$V\subset\C^2$ we need, firstly, to find an appropriate Faddeev type Green
function for $\bar\pa(\pa+\lambda dz_1)$ on $V$. One can check that for the
case $V=\C$ the Faddeev--Green function $g(z,\lambda)$ is a composition of
Cauchy--Green--Pompeiu kernels for the operators $f\mapsto\v=\bar\pa f$ and
$u\mapsto f=(\pa+\lambda dz)u$, where $u$, $f$, and $\v$ are respectively
a function, a (1,0)-form, and a (1,1)-form on $\C$. More precisely, one has the formula
$$g(z,\lambda)={-1\over i(2\pi)^2}\int\limits_{w\in\C}
{e^{\lambda w-\bar\lambda\bar w}dw\wedge d\bar w\over (w+z)\cdot\bar w}.$$
The main purpose of this paper is to construct an analogue of the Faddeev--Green
function on the Riemann surface $V$. To do this we need to find explicit
formulas $f=\hat R\v$ and $u=R_{\lambda}f$ (with appropriate estimates), for
solutions of the two equations $\bar\pa f=\v$ and
$(\pa+\lambda dz_1)u=f-{\cal H}_{\lambda}f$ on $V$. Here we consider $V$ equipped with
the Euclidean volume form  $d d^c|z|^2$, and we require
$\v\in L_{1,1}^1(V)$, $f\in {\tilde W}_{1,0}^{1,\tilde p}(V)$,
and $u\in L^{\infty}(V)$, $\tilde p>2$, with ${\cal H}_{\lambda}f$ being the projection
of $f$ on the subspace of pseudoholomorphic (1,0)-forms on $\tilde V$:
$\bar\pa {\cal H}_{\lambda}f=\bar\lambda d\bar z_1\wedge {\cal H}_{\lambda}f$.
\smallskip

The new formulas obtained in this paper for solution of $\bar\pa f=\v$ and
$(\pa+\lambda dz_1)u=f$ on $V$ one can interprete as explicit and more precise versions
of the classical Hodge--Riemann decomposition results on Riemann
surfaces. We
will define the Faddeev type Green function for $\bar\pa(\pa+\lambda dz_1)$
on $V$ as the kernel $g_{\lambda}(z,\xi)$ of the integral operator
$R_{\lambda}\circ \hat R$.
\bigskip

{\it Further results}
\smallskip

Let $\sigma\in C^{(2)}(V)$, with $\sigma>0$ on $V$, and
$\sigma\equiv{\rm const}$ on a
neighborhood of $\tilde V\b V$. Let $a_1,\ldots,a_g$ be generic points in this
neighborhood, with $g$ being the genus of $\tilde V$.
Using the Faddeev type Green function constructed here, we have in [HM]
obtained natural analogues of all steps of the Novikov reconstruction scheme
on the Riemann surface $V$.
In particular, under a smallness assumption on $d\log\sqrt{\sigma}$,
the existence (and uniqueness) of the solution $\mu(z,\lambda)$ of the Faddeev type
integral equation
$$\mu(z,\lambda)=1+{i\over 2}\int\limits_{\xi\in V}g_{\lambda}(z,\xi)
{\mu(\xi,\lambda)
d\,d^c\sqrt{\sigma}\over \sqrt{\sigma}}
+i\sum_{l=1}^gc_lg_{\lambda}(z,a_l),\ z\in V,\ \lambda\in\C$$
holds for any a priori fixed constants $c_1,\ldots,c_g$.
However (and this was overlooked in [HM]), there exists only one unique choice of
constants $c_l=c_l(\lambda,\sigma)$ for which the integral equation above
is equivalent to the differential equation
$$\bar\pa(\pa+\lambda dz_1)\mu={i\over 2}\bigl({d d^c\sqrt{\sigma}\over
\sqrt{\sigma}}\mu\bigr)+i\sum_{l=1}^gc_l\delta(z,a_l),$$
where $\delta(z,a_l)$ are Dirac measures concentrated in the points $a_l$
(see also $\S 4$ below).
\bigskip

\noindent {\bf $\S 1$. A Cauchy--Pompeiu type formula on an affine algebraic Riemann surface}
\smallskip

By $L_{p,q}(V)$ we denote the space of (p,q)-forms on $V$ with coefficients
in distributions of measure type on $V$. By
$L^s_{p,q}(V)$ we denote the space of (p,q)-forms on $V$ with absolutely
integrable in degree $s\ge 1$ coefficients with respect to the Euclidean
volume form on $V$. If $V=\C$ and $f$ is a function from $L^1(\C)$ such that
$\bar\pa f\in L_{0,1}(\C)$, then the generalized Cauchy formula has the
following form
$$f(z)=-{1\over 2\pi i}\int\limits_{\xi\in\C}{\bar\pa f(\xi)\wedge d\xi\over
{\xi-z}},\ \ z\in\C.$$
This formula becomes the classical Cauchy formula, when $f=0$ on $\C\b\Omega$
and $f\in {\cal O}(\Omega)$, where $\Omega$ is some bounded domain with
rectifiable boundary in $\C$.
The generalized Cauchy formula was discovered by Pompeiu [P1] in connection with his
solution of the Painlev\'e problem, i.e., in proving the existence for a totally
disconnected compact set $E$ with positive Lebesgue measure of a non-zero
function $f\in {\cal O}(\C\b E)\cap C(\C)\cap L^1(\C)$. The
Cauchy--Pompeiu formula has a large number of fundamental applications: in the theory of
distributions
(L.Schwartz), in approximations problems (E.Bishop, S.Mergelyan,
A.Vitushkin), in the solution of the corona problem (L.Carleson), in the theory of pseudo-analytic functions
(L.Bers, I.Vekua), and in inverse scattering and integrable equations (R.Beals,
R.Coifman, M.Ablowitz, D.Bar Yaacov, A.Fokas).
\smallskip

Motivated by applications to Electrical Impedance Tomography we develop in
this paper the Cauchy--Pompeiu type formulas on affine algebraic Riemann
surfaces $V\subset\C^2$ and give some applications.
\smallskip

Let $\tilde V$ be a smooth algebraic curve in ${\C}P^2$ given by the equation
$$\tilde V=\{w\in {\C}P^2;\ \tilde P(w)=0\},$$
with $\tilde P$ being a homogeneous holomorphic polynomial in the homogeneous
coordinates
$w=(w_0:w_1:w_2)\in {\C}P^2$. Without loss of generality
we may suppose that
\smallskip
\item{  i)} $\tilde V$ intersects ${\C}P^1_{\infty}=\{w\in {\C}P^2;\ w_0=0\}$
transversally, $\tilde V\cap {\C}P^1_{\infty}=\{a_1,\ldots,a_d\}$,
$d=\deg\tilde P$;
\smallskip
\item{ ii)} $V=\tilde V\b {\C}P^1_{\infty}$ is a connected curve in $\C^2$
with equation $V=\{z\in\C^2;\ P(z)=0\}$, where $P(z)=\tilde P(1,z_1,z_2)$
such that
$$\bigg|{\pa P/\pa z_1\over \pa P/\pa z_2}\bigg|\le {\rm const}(V),\ \ {\rm if}\ \
 |z_1|\ge r_0={\rm const}(V);$$
 \smallskip
\item{iii)} For any $z^*\in V$, such that ${\pa P\over \pa z_2}(z^*)=0$
we have ${\pa^2 P\over \pa z_2^2}(z^*)\ne 0$.
\smallskip

\noindent By the Hurwitz--Riemann theorem the number of such ramification points is equal to
$d(d-1)$. Let us equip $V$ with the Euclidean volume form $d\,d^c|z|^2$.
\bigskip

{\it Notation}
\smallskip

Let ${\tilde W}^{1,\tilde p}(V)=\{F\in L^{\infty}(V);
\bar\pa F\in L_{0,1}^{\tilde p}(V)\}$, $\tilde p>2$.
Let us denote by $H_{0,1}^p(V)$ the subspace in $L_{0,1}^p(V)$, $1<p<2$,
consisting of antiholomorphic forms. For all $p\in (1,2)$, the space
$H_{0,1}^p(V)$ coincides with the space of antiholomorphic forms on $V$
admitting an antiholomorphic extension to the compactification
$\tilde V\subset {\C}P^2$. Hence, by  the Riemann--Clebsch theorem
one has $\dim_{\C}H_{0,1}^p(V)=(d-1)(d-2)/2$ for all $p\in (1,2)$.
\bigskip

\noindent
{\bf Proposition 1.}
{\sl
Let $\{V_j\}$ be the connected components of
$\{z\in V;\ |z_1|>r_0\}$. Then for all $j\in \{1,\ldots,d\}$
there exist operators $R_1\colon\ L_{0,1}^p(V)\to L^{\tilde p}(V)$ and
$R_0\colon\ L_{0,1}^p(V)\to {\tilde W}^{1,\tilde p}(V)$ and
${\cal H}\colon\ L_{0,1}^p(V)\to H_{0,1}^p(V)$, $1<p<2$, $1/\tilde p=1/p-1/2$
such that, for all $\Phi\in L_{0,1}^p(V)$, one has the decomposition
$$\Phi=\bar\pa R\Phi+{\cal H}\Phi,\ \ {\rm where}\ \ R=R_1+R_0,
\eqno(1.1)$$
$$\eqalign{
&R_1\Phi={1\over 2\pi i}\int\limits_{\xi\in V}\Phi(\xi)
{d\xi_1\over {\pa P\over \pa\xi_2}}det\biggl[{\pa P\over \pa\xi}(\xi),
{{\bar\xi-\bar z}\over |\xi-z|^2}\biggr],\cr
&{\cal H}\Phi=\sum_{j=1}^g\bigl(\int\limits_V\Phi\wedge\omega_j\bigr)
\bar\omega_j,\cr}\eqno(1.2)$$
with $\{\omega_j\}$ being an orthonormal basis for the holomorphic (1,0)-forms on
$\tilde V$, i.e.,
$$\int\limits_V\omega_j\wedge\bar\omega_k=\delta_{jk},\ \ j,k=1,2,\ldots,g\,,$$
and
$$\lim_{\scriptstyle z\in V_j \atop \scriptstyle z\to\infty}R\Phi(z)=0.$$}
\bigskip

\noindent
{\it Remark 1.}
If $p\in [1,2)$ and $q\in (2,\infty]$ the condition
$\Phi\in L_{0,1}^p(V)\cap L_{0,1}^q(V)$ implies that
$R\Phi\in C(\tilde V)$.
\bigskip

\noindent
{\it Remark 2.}
For the case when $V=\C=\{z\in\C^2;\ z_2=0\}$ Proposition 1 and Remark 1
are reduced to the classical results of Pompeiu [P1], [P2] and of
Vekua [V].
\bigskip

\noindent
{\it Remark 3.}
Based on the technique of [HP] one can construct an explicit formula not only
for the main part $R_1$ of the $R$-operator, but for the whole operator $R$.
\bigskip

\noindent {\it Proof of Proposition 1:}
Let $Q(\xi,z)=\{Q_1(\xi,z),Q_2(\xi,z)\}$ be a pair of holomorphic
polynomials in the variables
$\xi=(\xi_1,\xi_2)$ and $z=(z_1,z_2)$, such that
$$\eqalign{
&Q(\xi,\xi)={\pa P\over \pa\xi}(\xi)\ \ {\rm and}\cr
&P(\xi)-P(z)=Q_1(\xi,z)(\xi_1-z_1)+Q_2(\xi,z)(\xi_2-z_2)\buildrel \rm def\over
 =\langle Q(\xi,z),\xi-z\rangle.\cr}\eqno(1.3)$$
The conditions i) and ii) imply that for $\ep_0$ small enough there exists
a holomorphic retraction $z\to r(z)$ of the domain
${\cal U}_{\ep_0}=\{z\in\C^2:\ |P(z)|<\ep_0\}$ onto the curve $V$.
\smallskip

Put ${\cal U}_{\ep,r}=\{z\in\C^2;\ |P(z)|<\ep,\ |z_1|<r\}$, where $0<\ep\le\ep_0$
and $r\ge r_0$. Put also $V^c=\{z\in\C^2;\ P(z)=c\}$, where $c\in\C$, $|c|\le\ep_0$ and
 $\tilde\Phi(z)=\Phi(r(z))$, $z\in {\cal U}_{\ep_0}$.
The condition $\Phi\in L_{0,1}^p(V)$ and properties of the retraction
$z\to r(z)$ together imply that $\bar\pa\tilde\Phi=0$ on ${\cal U}_{\ep_0}$ and
$$\|\tilde\Phi\|_{L^p(V^c)}\le {\rm const}(V)\cdot\|\Phi\|_{L^p(V)},\eqno(1.4)$$
uniformly with respect to $c$, for $|c|\le\ep_0$.
By results from [H] and [Po] we can choose the following explicit solution
$\tilde F_{\ep,r}$ on ${\cal U}_{\ep,r}$ of the $\bar\pa$-equation
$\bar\pa\tilde F_{\ep,r}=\tilde\Phi\big|_{{\cal U}_{\ep},r}$:

$$\eqalign{
&\tilde F_{\ep,r}(z)=\bigl({1\over 2\pi i}\bigr)
\bigg\{\int\limits_{\xi\in {\cal U}_{\ep,r}}
\tilde\Phi\wedge \det\bigg[{{\bar\xi-\bar z}\over |\xi-z|^2},
\bar\pa_{\xi}{{\bar\xi-\bar z}\over |\xi-z|^2}\bigg]
\wedge d\xi_1 \wedge d\xi_2+\cr
&\int\limits_{\xi\in b{\cal U}_{\ep,r}:\ |\xi_1|=r}
\tilde\Phi\wedge \bigg[-{(\bar\xi_2-\bar z_2)\over (\xi_1-z_1)|\xi-z|^2}\bigg]
\wedge d\xi_1 \wedge d\xi_2+\cr
&\int\limits_{\xi\in b{\cal U}_{\ep,r}:\ |P(\xi)|=\ep}
\tilde\Phi\wedge \det\bigg[{{\bar\xi-\bar z}\over |\xi-z|^2},
{Q\over {P(\xi)-P(z)}}\bigg]
\wedge d\xi_1 \wedge d\xi_2\bigg\},\ \ z\in {\cal U}_{\ep,R}.\cr}\eqno(1.5)$$
The property (1.4) implies that for any $z\in V$ we have
$$\eqalign{
&\int\limits_{\xi\in {\cal U}_{\ep,r}}
\tilde\Phi\wedge \det\bigg[{{\bar\xi-\bar z}\over |\xi-z|^2},
\bar\pa_{\xi}{{\bar\xi-\bar z}\over |\xi-z|^2}\bigg]
\wedge d\xi_1 \wedge d\xi_2\to 0,\ \ep\to 0\ \ {\rm and}\cr
&\int\limits_{\xi\in b{\cal U}_{\ep,r}:\ |\xi_1|=r}
\tilde\Phi\wedge \bigg[-{(\bar\xi_2-\bar z_2)\over (\xi_1-z_1)|\xi-z|^2}\bigg]
\wedge d\xi_1 \wedge d\xi_2\to 0,\ \ep\to 0,\ r\to\infty.\cr}$$
Hence for all $z\in V$ there exists
$\lim\limits_{\scriptstyle \ep\to 0 \atop\scriptstyle r\to\infty}
\tilde F_{\ep,r}=\tilde F(z)$, where
$$\tilde F(z)=-{1\over 2\pi i}\int\limits_{\xi\in V}
{\Phi d\xi_1\over {\pa P\over \pa\xi_2}(\xi)}\wedge \det
\bigg[{{\bar\xi-\bar z}\over |\xi-z|^2},Q(\xi,z)\bigg].\eqno(1.6)$$
From (1.5) and (1.6)   it follows that
$$\bar\pa_z\tilde F\big|_V=\Phi(z).\eqno(1.7)$$
Now put $F_1=R_1\Phi$. Using (1.2), (1.3),  (1.6), and (1.7) we obtain
$$\eqalign{
&\bar\pa_zF_1(z)\big|_V={1\over 2\pi i}\int\limits_{\xi\in V}
\Phi(\xi)\wedge {d\xi_1\over {\pa P\over \pa\xi_2}}\wedge
\det\bigg[{\pa P\over \pa\xi}(\xi),\bar\pa_z{{\bar\xi-\bar z}\over |\xi-z|^2}
\bigg]=\cr
&\Phi+{1\over 2\pi i}\int\limits_{\xi\in V}
\Phi(\xi)\wedge {d\xi_1\over {\pa P\over \pa\xi_2}}\wedge {1\over |\xi-z|^4}
\det\left|\matrix{
{\pa P\over \pa\xi_1}(\xi)&\ \xi_2-z_2\hfill\cr
{\pa P\over \pa\xi_2}(\xi)&-(\xi_1-z_1)\hfill\cr}\right|
\det\left|\matrix{
\bar\xi_1-\bar z_1 &d\bar z_1\hfill\cr
\bar\xi_2-\bar z_2 &d\bar z_2\hfill\cr}\right|=\cr
&\Phi+K\Phi,\cr}$$
where
$$K\Phi={1\over 2\pi i}\int\limits_{\xi\in V}
{\Phi(\xi)\wedge d\xi_1\over |\xi-z|^4}\wedge
{\langle{\pa P\over \pa\xi}(\xi),\xi-z\rangle\cdot
\langle{\pa\bar P\over \pa\bar\xi}(z),\bar\xi-\bar z\rangle\over
{\pa P\over \pa\xi_2}(\xi)\cdot {\pa\bar P\over \pa\bar\xi_2}(z)}
d\bar z_1.\eqno(1.8)$$
The estimate $R_1\Phi=F_1\in L^{\tilde p}(V)$ follows from the property
$\Phi\in L^p(V)$ and the following estimate of the kernel for the operator
$R_1$:
$$\bigg|\biggl({\pa P\over \pa\xi_2}(\xi)\biggr)^{-1}
\det\bigg[{\pa P\over \pa\xi}(\xi),{{\bar\xi-\bar z}\over |\xi-z|^2}
\bigg]d\xi_1\bigg|=O\biggl({1\over |\xi-z|}\biggr)
(|d\xi_1|+|d\xi_2|),$$
where $\xi,z\in V$.

For the kernels of operators $\Phi\mapsto K\Phi$ and $\Phi\mapsto \pa_zK\Phi$
we have the corresponding estimates
$$\eqalign{
&\bigg|{\langle{\pa P\over \pa\xi}(\xi),\xi-z\rangle\cdot
\langle{\pa\bar P\over \pa\bar\xi}(z),\bar\xi-\bar z\rangle \,d\xi_1\wedge d\bar z_1\over
{\pa P\over \pa\xi_2}(\xi)\cdot |\xi-z|^4\cdot
{\pa\bar P\over \pa\bar\xi_2}(z)}\bigg|=\cr \cr
&\left\{\matrix{
O\bigl({1\over {1+|z|^2}}\bigr)
|(d\xi_1+d\xi_2)\wedge (d\bar z_1+d\bar z_2)|\ \ &\ {\rm if}\ \ |\xi-z|\le 1,
\hfill\cr
O\bigl({1\over {|\xi-z|^2}}\bigr)
|(d\xi_1+d\xi_2)\wedge (d\bar z_1+d\bar z_2)|\ \ &\ {\rm if}\ \ |\xi-z|\ge 1.
\hfill\cr}\right.\cr}\eqno(1.9)$$
\smallskip
$$\eqalign{
&\bigg|\pa_z{\langle{\pa P\over \pa\xi}(\xi),\xi-z\rangle\cdot
\langle{\pa\bar P\over \pa\bar\xi}(z),\bar\xi-\bar z\rangle\,d\xi_1\wedge d\bar z_1\over
{\pa P\over \pa\xi_2}(\xi)\cdot |\xi-z|^4\cdot
{\pa\bar P\over \pa\bar\xi_2}(z)}\bigg|=\cr \cr
&\left\{\matrix{
O\bigl({1\over (1+|z|^2)|\xi-z|}\bigr)
|(d\xi_1+d\xi_2)\wedge (d\bar z_1+d\bar z_2)\wedge (dz_1+dz_2)|\ \
&{\rm if}\ \ |\xi-z|<1,\hfill\cr \cr
O\bigl({1\over |\xi-z|^3}\bigr)
|(d\xi_1+d\xi_2)\wedge (dz_1+dz_2)\wedge (d\bar z_1+d\bar z_2)|\ \
&{\rm if}\ \ |\xi-z|\ge 1,\hfill\cr}\right.\cr}\eqno(1.10)$$
 $\xi,z\in V$.
 \bigskip

These estimates imply that, for all $\tilde p>2$ and $p>1$, one has
$$\Phi_0\buildrel \rm def \over = K\Phi\in W_{0,1}^{1,\tilde p}(V)\cap
L_{0,1}^p(V).\eqno(1.11)$$
From estimates (1.9)-(1.11) it follows that the (0,1)-form $\Phi_0=K\Phi$ on $V$
can be considered also as a (0,1)-form on the compactification $\tilde V$ of $V$
in ${\C}P^2$ belonging to the spaces $W_{0,1}^{1,p}(\tilde V)$ for all $p<2$,
where $\tilde V$ is equipped with the projective volume form
$d\,d^c\ln(1+|z|^2)$.
\smallskip

From the Hodge--Riemann decomposition theorem [Ho], [W] we have
$$\Phi_0=\bar\pa(\bar\pa^*G\Phi_0)+{\cal H}\Phi_0,\eqno(1.12)$$
where ${\cal H}\Phi_0\in H_{0,1}(\tilde V)$, and $G$ is the Hodge--Green operator
for the Laplacian $\bar\pa\bar\pa^*+\bar\pa^*\bar\pa$ on $\tilde V$
with the properties
$$G(H_{0,1}(\tilde V))=0,\ \ \bar\pa G=G\bar\pa,\ \ \bar\pa^*G=G\bar\pa^*.$$
The decomposition (1.12) implies that
$$\bar\pa^*G\Phi_0\in W^{2,p}(\tilde V),\ p\in (1,2)\ \ {\rm and}\ \
{\cal H}\Phi_0\in H_{0,1}(\tilde V),$$
and this in turn implies that $\bar\pa^*G\Phi_0\in C(\tilde V)$.
Returning to the affine curve $V$ with the Euclidean volume form, we obtain that
$$\eqalign{
&\tilde R_0\Phi\buildrel \rm def \over =\bar\pa^*GK\Phi\big|_V\in
\tilde W^{1,\tilde p}(V),\ \forall\tilde p>2,\ \ {\rm where}\cr
&\tilde W^{1,\tilde p}(V)\buildrel \rm def \over =\{F\in L^{\infty}(V);\
\bar\pa F\in L^{\tilde p}_{0,1}(V)\},\cr
&{\rm and}\ \
{\cal H}\Phi\buildrel \rm def \over ={\cal H}K\Phi\big|_V\in
H_{0,1}^p(V),\ p>1.\cr}\eqno(1.13)$$

\noindent Now put $\tilde R=R_1+\tilde R_0$. Then, for all $\Phi\in L^p_{0,1}(V)$, we have
$\tilde R_0\Phi\in \tilde W^{1,\tilde p}(V)$, and
$\tilde R\Phi\in L^{\infty}(V)\cup L^{\tilde p}(V)$.
\smallskip

By Corollary 1.1 below, which is based only on (1.13), it follows that, for any form $\Phi\in L^p_{0,1}(V)$,
one has a limit

$$\lim\limits_{\scriptstyle z\to\infty \atop \scriptstyle z\in V_j}
\tilde R\Phi(z)\buildrel \rm def \over = {\tilde R}\Phi(\infty_j)\,.$$
Put $R_0\Phi=\tilde R_0\Phi-{\tilde R}\Phi(\infty_j)$ and
$R\Phi=\tilde R\Phi-{\tilde R}\Phi(\infty_j)$.
We then have property (1.1) for $R=R_1+R_0$.
This concludes the proof of Proposition 1.
\bigskip

\noindent
{\bf Corollary 1.1.}
{\sl
Let $F\in L^{\infty}(V)$ and $\bar\pa F\in L_{0,1}^p(V)$, $1<p<2$. Then,
for all $j\in \{1,\ldots,d\}$, there exists a limit
$\lim\limits_{\scriptstyle z\to\infty \atop\scriptstyle z\in V_j}F(z)
\buildrel \rm def \over = F(\infty_j)$
 such that  $(F-F(\infty_j))\big|_{V_j}\in L^{\tilde p}$.}
 \bigskip

\noindent
{\it Proof:}
Put $\bar\pa F=\Phi$. Then by  (1.13) we have
$\tilde R\Phi\in L^{\infty}(V)\cup L^{\tilde p}(V)$ and
$\bar\pa(F-\tilde R\Phi)={\cal H}\Phi$.
Then the function $h=F-\tilde R\Phi$ is harmonic on $V$. The estimates
$F\in L^{\infty}(V)$ and
$\tilde R\Phi\in L^{\tilde p}(V)\cup L^{\infty}(V)$
imply
by the Riemann extension theorem that $h$ can be extended to a
harmonic function
$\tilde h$ on $\tilde V$. Hence, $h=F-{\tilde R}\Phi\equiv {\rm const}=c$. This implies
that there exists $\lim\limits_{\scriptstyle z\to\infty \atop\scriptstyle z\in V_j}F(z)=
c_j\buildrel \rm def \over =F(\infty_j)$.
Corollary 1.1 is proved.
\bigskip

Corollary 1.1 admits the following useful reformulation.
\bigskip

\noindent
{\bf Corollary 1.2.}
{\sl
In the notations of Proposition 1, for any bounded function $\psi$ on $V$,
such that $\bar\pa\psi\in L^p(V)$, $1<p<2$, the following formula is valid:
$$\psi(z)=\psi(\infty_j)+R_0\bar\pa\psi+{1\over 2\pi i}\int\limits_{\xi\in V}
\bar\pa_{\xi}\biggl({\det\big[{\pa P\over \pa\xi}(\xi),\bar\xi-\bar z\big]d\xi
\over {\pa P\over \pa\xi_2}(\xi)\cdot  |\xi-z|^2}\biggr)\wedge\psi,$$
where $R_0\bar\pa\psi\in {\tilde W}^{1,\tilde p}(V)$, $1/\tilde p=1/2-1/p$,
and $R_0\bar\pa\psi(z)\to 0$, for $z\in V_j$, with $z\to\infty$.}
\bigskip

\noindent
{\bf $\S 2$. Kernels and estimates for $\bar\pa f=\v$ with
$\v\in L_{1,1}^1(V)$}
\smallskip

\noindent Let $\v$ be a (1,1)-form of class $L_{1,1}^{\infty}(V)$ with
support in $V_0=\{z\in V;\ |z_1|\le r_0\}$, where $r_0$ satisfies the condition
ii)  of $\S 1$.
\smallskip

If $V=\C$ then by classical results from [P1] and [V], the Cauchy--Pompeiu operator
$$\v\mapsto {dz\over 2\pi i}\int\limits_{\xi\in V_0}{\v(\xi)\over {\xi-z}}
\buildrel \rm def \over =\hat R\v$$
determines a solution $f=\hat R\v$ for the equation $\bar\pa f=\v$ on $\C$
with the property
$$f\in W_{1,0}^{1,\tilde p}(\C)\cap {\cal O}_{1,0}(\C\b V_0)\quad
{\rm for\ all} \quad \tilde p>2.$$
In this section we derive an analogous result for the case of an affine
algebraic Riemann surface $V\subset\C^2$.
\smallskip

Let $V\b V_0=\cup_{j=1}^dV_j$, where $\{V_j\}$ are the
connected components of $V\b V_0$.
\bigskip

\noindent
{\bf Lemma 2.1.}
{\sl
Let $\Phi=dz_1\rfloor\v$ and $f=Fdz_1=(R\Phi)dz_1$, where
$R$ is the operator from Proposition 1. Then
\smallskip
\item{  i)} $\Phi\in L_{0,1}^p(V_0)$, $p\in [1,2)$, $\Phi=0$ on $V\b V_0$;
\smallskip
\item{ ii)} $F\big|_{V_0}\in W^{1,p}(V_0)\ \forall p\in (1,2),\
f\in {\tilde W}_{1,0}^{1,\tilde p}(V)\ \forall\tilde p\in (2,\infty)$,
\item{    } $\bar\pa F=\Phi-{\cal H}\Phi$, where ${\cal H}$ is the operator
from Proposition 1,
\item{    } $\bar\pa f=\v-({\cal H}\Phi)\wedge dz_1$, and
\item{    } $\|F\|_{L^{\infty}(V\b V_0)}+
\|F\|_{L^{\tilde p}(V_0)}\le {\rm const}(V,p)\|\Phi\|_{L^p(V)},\ \
1/\tilde p=1/p-1/2$;
\smallskip

\item{iii)} if, in addition, $\v\in W^{1,\infty}(V)$, then
$f\in {\tilde W}^{2,\tilde p}(V)$.}
\bigskip

\noindent
{\it Proof:}
\smallskip

\item{  i)} The property $\Phi\big|_{V\b V_0}=0$ follows from
$\v\big|_{V\b V_0}=0$.

Put $V_0^{\pm}=\{z\in V_0;\ \pm\big|{\pa P\over \pa z_2}\big|\ge \pm
\big|{\pa P\over \pa z_1}\big|\}$.
The definition $\Phi=dz\rfloor\v$ implies that
$$\eqalign{
&\Phi\big|_{V_j^+}=\Phi^+d\bar z_1,\ \ {\rm where}\ \
\Phi^+\in L^{\infty}(V_0^+)\ \ {\rm and}\cr
&\Phi\big|_{V_j^-}=\Phi^-d\bar z_2/(\pa z_1/\pa z_2),\ \  {\rm where}\ \
\Phi^-\in L^{\infty}(V_0^-).\cr}\eqno(2.1)$$
\item{}The properties (2.1) imply that $\Phi\in L_{0,1}^p(V_0)$ for all $p\in (1,2)$.
\smallskip
\item{ ii)} The equalities $\bar\pa F=\Phi-{\cal H}\Phi$ and
$\bar\pa f=\v-({\cal H}\Phi)\wedge dz_1$ follow from Proposition~1 together with
the definitions $\Phi=dz_1\rfloor\v$ and $f=Fdz_1$. The inclusions
$F\in L^{\infty}(V\b V_0)$ and $F\big|_{V_0}\in W^{1,p}(V_0)$ follow from
the formula $F=R\Phi$ and Proposition 1. The inclusion
$f\in {\tilde W}^{1,\tilde p}_{1,0}(V)$ follows from the equalities
$\bar\pa f=\v-({\cal H}\Phi)\wedge dz_1$, $f=Fdz_1$ and Proposition 1.
\smallskip
\item{iii)} If, in addition, $\v\in W^{1,\infty}_{1,1}(V)$, then
$f\in {\tilde W}^{2,\tilde p}_{1,0}(V)$. It follows from equalities above
with $\v\in W^{1,\infty}_{1,1}(V)$ and ${\rm supp}\,\v\subset V_0$.
\bigskip

\noindent
{\bf Lemma 2.2.}
{\sl
For each $(0,1)$-form $g\in H_{0,1}^p(V)$ there exists a $(1,0)$-form $h\in L_{1,0}^p(\tilde V)$ ($1\le p<2$),
unique up to adding
holomorphic (1,0)-forms on $\tilde V$, such that
$$\bar\pa h\big|_{\tilde V}=gdz_1.\eqno(2.2)$$}

\noindent
{\it Proof:}
For any $g\in H_{0,1}^p(V)$ the (1,1)-form $g\wedge dz_1$ determines a
current $G$ on $\tilde Y$ by the equality
$$\langle G,\chi\rangle\buildrel \rm def \over =\lim\limits_{R\to\infty}
\sum\limits_{j=1}^d\big[\int\limits_{V_j}(\chi-\chi_j(\infty))gdz_1+
\chi_j(\infty)\int\limits_{\{z\in V_j;\ |z_1|<r\}}g\wedge dz_1\big],$$
where $\chi\in C^{(\ep)}(\tilde V)$, $\ep>0$ and $\chi_j(\infty)=
\lim\limits_{\scriptstyle z\in V_j \atop \scriptstyle z\to\infty}\chi(z)$.
\smallskip

\noindent By Serre duality [S], the current $G$ is $\bar\pa$-exact on $\tilde V$ if and only if
$$\langle G,1\rangle=\lim\limits_{R\to\infty}\int\limits_{\{z\in V;\ |z_1|\le r\}}
g\wedge dz_1=0.\eqno(2.3)$$
Let us check (2.3). We have
$$\int\limits_{\{z\in V:\ |z_1|\le r\}}g\wedge dz_1=-
\int\limits_{\{z\in V:\ |z_1|=r\}}z_1\wedge g.$$
Putting $w_1=1/z_1$ into the right-hand side of this equality, we obtain
$$\int\limits_{\{z\in V:\ |z_1|\le r\}}g\wedge dz_1=
-\sum\limits_{j=1}^d\int\limits_{|w_1|=1/r}g_j(\bar w_1){d\bar w_1\over w_1}=0.
$$
Here the last equality follows from the properties
$$g_j(\bar w_1)d\bar w_1=g\big|_{V_j\cap\{|w_1|\le 1/r\}}\ \ {\rm and}\ \
\bar g_j\in {\cal O}(D(0,1/r)).$$
Hence, by (2.3) there exists $h\in L_{1,0}^1(\tilde V)$ such that equality
(2.2) is valid in the sense of currents. Moreover, any solution of (2.2)
automatically belongs to $L_{1,0}^p(\tilde V)$, $1<p<2$. Such a solution $h$
of (2.2) is unique up to holomorphic (1,0)-forms on $\tilde V$ because the
conditions $h\in L_{1,0}^p(\tilde V)$ and $\bar\pa h=0$ on $V$
imply that $h$ extends as a holomorphic (1,0)-form on $\tilde V$.
\bigskip

\noindent
{\it Notation:}
Let ${\cal H}^{\perp}\colon H_{0,1}^p(V)\to L_{1,0}^p(\tilde V)$ ($1<p<2$)
be the operator
defined by the formula $g\mapsto {\cal H}^{\perp}g$, where ${\cal H}^{\perp}g$
is the unique solution $h$ of (2.2) in $L_{1,0}^p(\tilde V)$ with the
property
$$\int\limits_Vh\wedge\tilde g=0\quad {\rm for\ all}\quad \tilde g\in H_{0,1}^p(V).$$
Lemma 2.2 guarantees the existence and uniqueness of
$H^{\perp}g\in L_{1,0}^p(\tilde V)$ for any $g\in H_{0,1}^p(V)$.
\bigskip

\noindent
{\bf Proposition  2.}
{\sl
Let $R$ be the operator defined by formula (1.1), and ${\cal H}$ the operator
defined by formula (1.13). For any (1,1)-form
$\v\in L_{1,1}^1(V)\cap L_{1,1}^{\infty}(V)$ with support in $V_0$, put
$$\hat R\v=R^1\v+R^0\v,\eqno(2.4)$$
where
$$R^1\v=(R(dz_1\rfloor\v)) dz_1,\ \
R^0\v={\cal H}^{\perp}\circ {\cal H}(dz_1\rfloor\v).$$
Then
$$\bar\pa\hat R\v=\v,\eqno(2.5),$$
$$\eqalign{
&f=Fdz_1=\hat R\v\in {\tilde W}^{1,\tilde p}_{1,0}(V)\ {\rm for\ all}\ \tilde p\in (2,\infty),
\ \ F\big|_{V_0}\in W^{1,p}(V_0)\ {\rm for\ all}\ p\in (1,2)\cr
&{\rm and}\ \
f\big|_{V_l}=\sum\limits_{k=1}^{\infty}{c_k^{(l)}\over z_1^k}dz_1+
b_ldz_1,\ \ {\rm if}
\ \ |z_1|\ge r_0.\cr}\eqno(2.6)$$
Here $ l=1,\ldots,d$, and $b_l=0$ for $l=j$.
}
\bigskip

\noindent
{\it Proof:}
The properties (2.5)   and  $f=\hat R\v\in {\tilde W}^{1,\tilde p}_{1,0}(V)$
follow from
Proposition 1 and Lemmas 2.1, 2.2. The properties (2.5) and
$\v\big|_{V\b V_0}=0$
imply analyticity of $f$ on $V\b V_0$. The series expansion (2.6) follows from
the analyticity of $f\big|_{V\b V_0}$ and the inclusion
$f\big|_{V\b V_0}\in L^{\infty}_{1,0}(V\b V_0)$.
\bigskip

\noindent
{\bf Supplement:}
Let $\tilde V_0=\{z\in V:\ |z_1|\le \tilde r_0\}$, where $\tilde r_0>r_0$.
If
${\rm supp}\,\v\subseteq V_0$ and
$$(\v-\sum_{l=1}^gc_l\delta(z,a_l))\in L^{\infty}_{1,1}(V),\ \
{\rm where}\ \ a_l\in V_{j(l)}\cap\tilde V_0,$$
then
$$(\hat R\v-\sum_{l=1}^gc_l\hat R(\delta(z,a_l)))
\in W^{1,\tilde p}_{1,0}(V).$$
\bigskip

\noindent
{\bf $\S 3$. Kernels and estimates for $(\pa+\lambda dz_1)u=f$, with
$f\in W_{1,0}^{1,\tilde p}(V)$}
\smallskip

If $V=\C$ then the equation $\pa u+\lambda udz_1=f$ was also introduced by
Pompeiu [P2]. One can check that this equation can be solved by the
explicit formula:
$$
e^{\lambda z-\bar\lambda\bar z}u(z)={1\over 2\pi i}\int\limits_{\xi\in\C}
{e^{\lambda\xi-\bar\lambda\bar\xi}f(\xi)d\bar\xi\over {\bar\xi-\bar z}}\
\buildrel \rm def \over = \ \lim\limits_{r\to\infty}{1\over 2\pi i}
\int\limits_{\{\xi\in\C:\ |\xi|<r\}}
{e^{\lambda\xi-\bar\lambda\bar\xi}f(\xi)d\bar\xi\over {\bar\xi-\bar z}}.$$

For a Riemann surface $V=\{z\in\C^2:\ P(z)=0\}$ we will obtain the
following generalization of this formula.
\bigskip

\noindent
{\bf Proposition 3.}
{\sl
Let $f=Fdz_1$ be a (1,0)-form   as in Proposition 2, i.e.,
$F\big|_{V_0}\in W^{1,p}(V_0)$ for all $p\in (1,2)$,
$f\in {\tilde W}^{1,\tilde p}_{1,0}(V)$ for all $\tilde p\in (2,\infty)$, and
${\rm supp}\,\bar\pa f\subset V_0$.
Let $e_{\lambda}(\xi)=e^{\lambda\xi_1-\bar\lambda\bar\xi_1}$.
Put
$$\overline{R_1(\bar e_{\lambda}\bar f)}\buildrel \rm def \over =
-{1\over 2\pi i}
\lim\limits_{r\to\infty}\int\limits_{\{\xi\in V:\ |\xi|<r\}}
e_{\lambda}(\xi)f(\xi){d\bar\xi_1\over \pa\bar P/\pa\bar\xi_2}
\det\big[{\pa\bar P\over \pa\bar\xi}(\xi),{{\xi-z}\over |\xi-z|^2}\big].$$
Put also ${\cal H}f\buildrel \rm def \over =\overline{{\cal H}\bar f}$,
where ${\cal H}$ is the operator from Proposition 1.
Finally, let
$$u=R_{\lambda}f=R_{\lambda}^1f+R_{\lambda}^0f,\eqno(3.1)$$
where
$R_{\lambda}^1f+R^0_{\lambda}f=
e_{-\lambda}(z)\cdot\overline{R_1(\bar e_{\lambda}\bar f)}+
e_{-\lambda}(z)\cdot\overline{R_0(\bar e_{\lambda}\bar f)}$, with
$R_1$ and $R_0$ being the operators from
Proposition 1.
\smallskip

Then for all $\lambda\ne 0$ one has:
\smallskip

\item{  i)} $(\pa+\lambda dz_1)R_{\lambda}f=f-{\cal H}_{\lambda}(f)$,
where
${\cal H}_{\lambda}(f)=e_{-\lambda}(z){\cal H}(e_{\lambda}f)$.
\smallskip

\item{ ii)}
$$\eqalign{
&\|u-u(\infty_l)\|_{L^{\infty}(V)}\le\cr
& {\rm const}(V,p)\cdot\min
\bigl({1\over \sqrt{|\lambda|}},
{1\over |\lambda|}\bigr)
\bigl(\|F\|_{L^{\tilde p}(V_0)}+\|F\|_{L^{\infty}(V\b V_0)}+
\|\pa F\|_{L^p_{1,0}(V)}\bigr),\cr
&\|\pa u\|_{L^{\tilde p}_{1,0}(V)}\le {\rm const}(V,p)\cdot
\|\pa F\|_{L^p_{1,0}(V)},\ \ {\rm where}\ \ 1/\tilde p=1/p-1/2,\
l=1,\ldots,d.\cr}$$

\item{iii)}
$$\eqalign{
&\|(1+|z_1|)(u-u(\infty_l))\|_{L^{\infty}(V_l)}\le
{{\rm const}(V,\tilde p)\over \sqrt{|\lambda|}}
(\|F\|_{L^{\tilde p}(V_0)}+\|F\|_{L^{\infty}(V\b V_0)})
,\cr
&\|(1+|z_1|)\pa u\|_{L^{\infty}_{1,0}(V)}\le
{\rm const}(V,\tilde p)(\sqrt{|\lambda|}+1)
(\|F\|_{L^{\tilde p}(V_0)}+\|F\|_{L^{\infty}(V\b V_0)})
,\forall\tilde p>2.\cr}$$
}
\bigskip

\noindent{\bf Supplement:}
Put
$$L^{2\pm\ep}(V)=\{u;\, u\big|_{\tilde V_0}\in L^{2-\ep}(\tilde V_0),\ \
u\big|_{V\b \tilde V_0}\in L^{2+\ep}(V\b \tilde V_0)\}.$$
If $f=f_0-f_1$, where $f_0\in {\tilde W}^{1,\tilde p}_{1,0}(V)$,
${\rm supp}\,\bar\pa f_0\subset V_0$ and $f_1=
\sum\limits_{l=1}^gc_l\hat R(\delta(z,a_l))$, $a_l\in V_{j(l)}\cap\tilde V_0$,
then instead of i)-iii) we have i) and the following conclusion:
\smallskip
{\sl \item{ii)$^{\prime}$}}
$$\eqalign{
&\|R_{\lambda}f-R_{\lambda}f(\infty_l)\|_{L^{2+\ep}(V_l)}\le
{{\rm const}(V,\tilde p)\over \ep}
\min\bigl(|\lambda|^{-1/2},|\lambda|^{-1}\bigr)
\bigl(\|f_0\|_{{\tilde W}^{1,\tilde p}_{1,0}}+\sum_{l=1}^g|c_l|\bigr),\cr
&\|\pa R_{\lambda}f\|_{L^{2\pm\ep}_{1,0}(V)}\le {{\rm const}(V,\tilde p)\over \ep}
\bigl(\|f_0\|_{{\tilde W}^{1,\tilde p}_{1,0}}+\sum_{l=1}^g|c_l|\bigr),\cr
&\|{\cal H}_{\lambda}(f)\|_{L^{\infty}_{1,0}(V)}\le {{\rm const}(V,\tilde p)\over
(1+|\lambda|)}
\bigl(\|f_0\|_{{\tilde W}^{1,\tilde p}_{1,0}}+\sum_{l=1}^g|c_l|\bigr), {\rm where} \ \tilde p>2, 0<\ep<1/2.\cr}$$
\bigskip

\noindent
{\bf Lemma 3.1.}
{\sl
Put
$$J(z)=\int\limits_{\{\xi\in\C:\ |\xi|<\rho\}}
{\psi(\xi)d\xi\wedge d\bar\xi\over |\xi|\cdot |\xi-z|},\ \ z\in\C,$$
where $\psi\in L^p(V_0)$, $p>1$. Then, for any $\ep>0$ and any $\tilde p>2$,
one has the estimate
$$\|J(z)\|_{L^{\tilde p}(\C)}\le
{1\over \ep}O(\rho^{(2-2\ep\tilde p)/\tilde p})\cdot
\|\psi\|_{L^{(1+\ep)/\ep}(V_0)}.$$
}
\bigskip

\noindent
{\it Proof:}
Using the notation $\|\psi\|_{\ep}=\|\psi\|_{L^{(1+\ep)/\ep}(V_0)}$, we obtain from
the expression for $J(z)$ the following estimates:
$$\eqalign{
&|J(z)|\le \biggl(\int\limits_{|\xi|\le\rho}{|d\xi\wedge d\bar\xi|\over
|\xi|^{1+\ep}|\xi-z|^{1+\ep}}\biggr)^{1/(1+\ep)}\cdot \|\psi\|_{\ep}\le\cr
&O\biggl(\int\limits_{r=0}^{\rho}{dr\over r^{\ep}}
\int\limits_0^1{d\v\over (|r-|z||+|z|\v)^{1+\ep}}\biggr)^{1/(1+\ep)}\cdot
\|\psi\|_{\ep}\le\cr
&O\biggl(\int\limits_{r=0}^{\rho}{dr\over r^{\ep}}
{1\over |z|^{1+\ep}}
\int\limits_0^1{d\v\over (|{r\over |z|}-1|+\v)^{1+\ep}}\biggr)^{1/(1+\ep)}
\cdot \|\psi\|_{\ep}\le\cr
&{1\over \ep}O\biggl(\int\limits_0^{\rho}{dr\over r^{\ep}}
{1\over |z|}
\biggl({1\over |r-|z||^{\ep}}-{1\over (|r-|z||+|z|)^{\ep}}
\biggr)^{1/(1+\ep)}\cdot \|\psi\|_{\ep}.\cr}$$
From the last estimate we deduce
$$\eqalign{
&|J(z)|\le {1\over \ep}O\biggl({1\over |z|}\biggl(
\int\limits_0^{|z|}{dr\over r^{\ep}|z|^{\ep}}+\int\limits_{|z|}^{\rho}
{\ep|z|\over r^{\ep}r^{\ep}r}\biggr)\biggr)^{1/(1+\ep)}\|\psi\|_{\ep},\ \
{\rm if}\ \ |z|\le\rho,\cr
&{\rm and}\ \ |J(z)|\le {1\over \ep}O\biggl({1\over |z|}
\int\limits_0^{\rho}
{dr\over r^{\ep}|z|^{\ep}}\biggr)\biggr)^{1/(1+\ep)}\|\psi\|_{\ep},\ \
{\rm if}\ \ |z|\ge\rho.\cr}$$
These equalities imply
$$\eqalign{
&|J(z)|\le {1\over \ep}O\biggl(\bigl({1\over |z|}\bigr)^{2\ep/(1+\ep)}\biggr)
\|\psi\|_{\ep},\ \ {\rm if}\ \ |z|\le\rho,\cr
&|J(z)|\le {1\over \ep}O\biggl(\bigl({1\over |z|}\rho^{(1-\ep)/(1+\ep)}\biggr)
\|\psi\|_{\ep},\ \ {\rm if}\ \ |z|\ge\rho.\cr}$$
Putting $|z|=t$ we obtain finally that
$$\|J\|_{L^{\tilde p}(\C)}\le {1\over \ep}O\biggl(\int\limits_0^{\rho}
{dt\over t^{2\ep\tilde p/(1+\ep)-1}}+\rho^{{{1-\ep}\over {1+\ep}}\tilde p}
\int\limits_{\rho}^{\infty}{dt\over t^{\tilde p-1}}\biggr)^{1/\tilde p}
\|\psi\|_{\ep}\le {1\over \ep}
O\biggl(\rho^{{{2-\ep\tilde p}\over \tilde p}}\biggr)\|\psi\|_{\ep}.$$
Lemma 3.1 is proved.
\bigskip

\noindent
{\it Proof of Proposition 3:}

\item{  i)}
$$\eqalign{
&(\pa+\lambda dz_1)R_{\lambda}f=(\pa+\lambda dz_1)
e_{-\lambda}(z)\cdot\overline{R(\bar e_{\lambda}\bar f)}=\cr
&\pa(e_{-\lambda}(z))\cdot\overline{R(\bar e_{\lambda}\bar f)}+
e_{-\lambda}(z)\pa(\overline{R(\bar e_{\lambda}\bar f)})+
\lambda dz_1e_{-\lambda}(z)\cdot\overline{R(\bar e_{\lambda}\bar f)}=\cr
&(-\lambda dz_1+\lambda dz_1)
e_{-\lambda}(z)\cdot\overline{R(\bar e_{\lambda}\bar f)}+\cr
&e_{-\lambda}(z)\cdot
(e_{\lambda}(z)f-\overline{{\cal H}{\bar e}_{\lambda}\bar f})=
f-e_{-\lambda}{\cal H}(e_{\lambda}f)\buildrel \rm def \over =
f-{\cal H}_{\lambda}f,\cr}$$
where we have used the equality (1.1) from Proposition 1.
\smallskip

\item{iii)} Let $r\ge r_0$. Let the functions $\chi_{\pm}\in C^{(1)}(V)$ be such that
$\chi_++\chi_-\equiv 1$ on $V$,
${\rm supp}\,\chi_+\subset\{\xi\in V: |\xi_1|<2r\}$,
${\rm supp}\,\chi_-\subset\{\xi\in V: |\xi_1|\ge r\}$, and $|d\chi_{\pm}|=O(1/r)$.
We then have $u=u_++u_-$, where
$$u_{\pm}(z)=R_{\lambda}(\chi_{\pm}f).\eqno(3.1)_{\pm}$$
Using the properties $f\in L^{\infty}(V)$ and $|e_{\lambda}|\equiv 1$, in combination with the
equality $\pa u_+=\chi_+Fdz_1-\lambda u_+dz_1-{\cal H}_{\lambda}(\chi_+f)$,
we obtain for $u_+$ and
${\pa u_{\pm}\over \pa z_1}$ the estimates:
$$\eqalign{
&\|(1+|z|)(u_+(z)-u_+(\infty_l))\|_{L^{\infty}(V_l)}=O(r)
\|f\|_{L^{\infty}_{1,0}(V)},\ l=1,\ldots,d,\cr
&\|(1+|z|)\pa u_+(z)\|_{L^{\infty}_{1,0}(V)}=
O(\lambda r+1)\|f\|_{L^{\infty}_{1,0}(V)}.\cr}\eqno(3.2)$$
In order to estimate $u_-$ we transform the expression $(3.1)_-$ using
the series expansion (2.6) for $f\big|_{V_j}$, and we  integrate by part. We
thus obtain
$$\eqalign{
&u_-(z)=R_{\lambda}\chi_-f=R_{\lambda}^1\chi_-f+R_{\lambda}^0\chi_-f=\cr
&-{e_{-\lambda}(z)\over 2\pi i}{1\over \lambda}\int\limits_{\xi\in V}
{e^{\lambda\xi_1-\bar\lambda\bar\xi_1}(d\chi_-)F\wedge d\bar\xi_1
\det\big[{\pa\bar P\over \pa\bar\xi}(\xi),\xi-z\big]
\over {\pa \bar P\over \pa\bar\xi_2}(\xi)\cdot  |\xi-z|^2}+\cr
&{e_{-\lambda}(z)\over 2\pi i}{1\over \lambda}\sum\limits_j
\int\limits_{\xi\in V_j}
e^{\lambda\xi_1-\bar\lambda\bar\xi_1}\chi_-
\bigl(\sum\limits_{k=1}^{\infty}k{c_k^{(j)}\over \xi_1^{k+1}}\bigr)
{d\xi_1\wedge d\bar\xi_1
\det\big[{\pa\bar P\over \pa\bar\xi}(\xi),\xi-z\big]
\over {\pa \bar P\over \pa\bar\xi_2}(\xi)\cdot  |\xi-z|^2}-\cr
&{e_{-\lambda}(z)\over 2\pi i}{1\over \lambda}
\int\limits_{\xi\in V}
e^{\lambda\xi_1-\bar\lambda\bar\xi_1}\chi_-
F\pa_{\xi}\biggl(
{\det\big[{\pa\bar P\over \pa\bar\xi}(\xi),\xi-z\big]d\bar\xi_1
\over {\pa \bar P\over \pa\bar\xi_2}(\xi)\cdot  |\xi-z|^2}\biggr)+
e_{-\lambda}(z)\bar R_0(e_{\lambda}\chi_-f),\cr}\eqno(3.3)$$

\noindent where the operator $R_0=\bar\pa^*GK$ is defined by  (1.13).
Using Corollary 1.2 we have, in addition,
$$\eqalign{
&-{e_{-\lambda}(z)\over 2\pi i}{1\over \lambda}
\int\limits_{\xi\in V}
e^{\lambda\xi_1-\bar\lambda\bar\xi_1}\chi_-
F\pa_{\xi}\biggl(
{\det\big[{\pa\bar P\over \pa\bar\xi}(\xi),\xi-z\big]d\bar\xi_1
\over {\pa \bar P\over \pa\bar\xi_2}(\xi)\cdot  |\xi-z|^2}\biggr)=\cr
&{e_{-\lambda}(z)\over 2\pi i}{1\over \lambda}e_{\lambda}(z)\chi_-(z)F(z)-
{e_{-\lambda}(z)\over 2\pi i}{1\over \lambda}
\bar R_0(\pa(e_{\lambda}\chi_-F))=\cr
&{1\over 2\pi i}{1\over \lambda}\chi_-(z)F(z)-
{e_{-\lambda}(z)\over 2\pi i}{1\over \lambda}
\bar R_0(\pa(e_{\lambda}\chi_-F)).\cr}$$
Putting the last equality in (3.3) and making use of the properties
$|e_{\lambda}|\equiv 1$, $|d\chi_-|=O(1/r)$,
$\pa u_-=\chi_-Fdz_1-\lambda u_-dz_1-{\cal H}_{\lambda}(\chi_-f)$, and the
property of $R_0$, we obtain
from Proposition,1:
$$\eqalign{
&\|(1+|z_1|)(u_--u_-(\infty_l))\|_{L^{\infty}(V_l)}=\cr
&O\bigl({1\over |\lambda|r}\bigr)
(\|F\|_{L^{\tilde p}(V_0)}+\|F\|_{L^{\infty}(V\b V_0)})
+\|(1+|z_1|)\bar R_0(e_{\lambda}\chi_-f)\|_{L^{\infty}(V)}+\cr
&{1\over 2\pi |\lambda|}
\|(1+|z_1|)\bar R_0\pa(e_{\lambda}\chi_-F)\|_{L^{\infty}(V)}\le
O\bigl({1\over |\lambda|r}\bigr)\|F\|_{L^{\tilde p}(V)},\ l=1,\ldots,d\cr
&{\rm and}\ \ \|(1+|z_1|){\pa u_-\over \pa z_1}\|_{L^{\infty}(V)}=
O(1/r+1)(\|F\|_{L^{\tilde p}(V_0)}+\|F\|_{L^{\infty}(V\b V_0)})
.\cr}\eqno(3.4)$$
The estimates (3.2) and (3.4) imply
$$\eqalign{
&\|(1+|z_1|)(u-u(\infty_l))\|_{L^{\infty}(V_l)}=\cr
&O\bigl(r+{1\over |\lambda|r}\bigr)
(\|F\|_{L^{\tilde p}(V_0)}+\|F\|_{L^{\infty}(V\b V_0)}),\cr
&{\rm and}\ \ \|(1+|z_1|)\pa u\|_{L^{\infty}_{1,0}(V)}=\cr
&O(|\lambda|r+1/r+1)(\|F\|_{L^{\tilde p}(V_0)}+
\|F\|_{L^{\infty}(V\b V_0)}),\ \forall\tilde p>2.\cr}\eqno(3.5)$$

\noindent Putting in (3.5) $r=r_0/\sqrt{|\lambda|}$ we obtain iii).
\smallskip

\item{ ii)} For proving ii) let us put $r=\tilde r_0$ and transform
$(3.1)_+$ for $u_+$ in the following way:
$$\eqalign{
&u_+(z)=R_{\lambda}\chi_+f=\cr
&-{e_{-\lambda}(z)\over 2\pi i}{1\over \lambda}\int\limits_{|\xi_1|\le r}
{e^{\lambda\xi_1-\bar\lambda\bar\xi_1}d\chi_+F\wedge d\bar\xi_1
det\big[{\pa\bar P\over \pa\bar\xi}(\xi),\xi-z\big]
\over {\pa \bar P\over \pa\bar\xi_2}(\xi)\cdot  |\xi-z|^2}-\cr
&{e_{-\lambda}(z)\over 2\pi i}{1\over \lambda}
\int\limits_{|\xi_1|\le r}
{e^{\lambda\xi_1-\bar\lambda\bar\xi_1}\chi_+
\pa F\wedge d\bar\xi_1
det\big[{\pa\bar P\over \pa\bar\xi}(\xi),\xi-z\big]
\over {\pa \bar P\over \pa\bar\xi_2}(\xi)\cdot  |\xi-z|^2}-\cr
&{e_{-\lambda}(z)\over 2\pi i}{1\over \lambda}
\int\limits_{|\xi_1|\le r}
e^{\lambda\xi_1-\bar\lambda\bar\xi_1}\chi_+
F\pa\biggl(
{det\big[{\pa\bar P\over \pa\bar\xi}(\xi),\xi-z\big]d\bar\xi_1
\over {\pa \bar P\over \pa\bar\xi_2}(\xi)\cdot  |\xi-z|^2}\biggr)+
e_{-\lambda}(z)\bar R_0(e_{\lambda}\chi_+f),\cr}\eqno(3.6)$$

\noindent where $R_0$ is the operator from Proposition 1.
Using the last expression for $u_+(z)$, together with the property
$F\big|_{V_0}\in W^{1,p}(V_0)$ and Corollary 1.2, we obtain
$$\|u_+\|_{L^{\infty}(V)}=O(1/\lambda)
\|F\|_{{\tilde W}^{1,p}(V_0)}.\eqno(3.7)$$
This inequality together with (3.4) and statement iii) proves the first
part of statement ii). Formula $u=R_{\lambda}f$ implies
$\pa_zu=f-\lambda dz_1u-{\cal H}_{\lambda}f$. From this and from the already
obtained estimates for $u$ we deduce the second part of statement ii):
$$\|\pa u\|_{L^{\tilde p}_{1,0}(V)}\le {\rm const}(V,p)
\|\pa F\|_{L^p_{1,0}(V)}.$$

\item{ii)$^{\prime}$} In order to prove in this case the estimate for
$u=R_{\lambda}f$ with $|\lambda|\le 1$, we combine the

\noindent arguments above with
Lemma 3.1, and obtain instead of  (3.5) the following:
$$\eqalign{
&\|u-u(\infty_l)\|_{L^{2+\ep}(V_l)}\le {1\over \ep}O\bigl(r+{1\over |\lambda| r}\bigr)
\bigl(\|f_0\|_{{\tilde W}^{1,\tilde p}_{1,0}(V)}+\sum_{l=1}^g|c_l|\bigr)\cr
&\|\pa u\|_{L^{2\pm\ep}_{1,0}(V)}\le {\lambda\over \ep}
O\bigl(r+{{1+r}\over |\lambda| r}\bigr)
\bigl(\|f_0\|_{{\tilde W}^{1,\tilde p}_{1,0}(V)}+\sum_{l=1}^g|c_l|\bigr).\cr}
\eqno(3.5)^{\prime}$$
Putting in (3.5)$^{\prime}$ $r=r_0\big/\sqrt{|\lambda|}$, we obtain
the required estimate for $R_{\lambda}f$ with $|\lambda|\le 1$.
To prove the estimate for $u=R_{\lambda}f$ with $|\lambda|\ge 1$, we use
 (3.6) and the Calderon--Zygmund $L^{2-\ep}$-estimate for the singular integral on
the right hand side of (3.6).
\smallskip

In order to prove the statement concerning ${\cal}H_{\lambda}f$, we just perform an
integration by parts in the expression
\vskip-.5cm
$${\cal H}_{\lambda}f=e_{-\lambda}{\cal H}(e_{\lambda}f)=
\sum_{l=1}^ge_{-\lambda}(z)
\biggl(\int\limits_{\tilde V}e_{\lambda}(\xi)f(\xi)\wedge
\overline{\omega_l(\xi)}\biggr)\omega_l(z),$$
\vskip-.5cm
\noindent where $f=f_0+\sum\limits_{l=1}^gc_l\hat R(\delta(z,a_l))$, and where
$\{\omega_l,l=1,\ldots, g\}$ is an orthonormal basis of holomorphic
(1,0)-forms on $\tilde V$.
\bigskip

\noindent
{\bf $\S 4$. Faddeev type Green function for $\bar\pa(\pa+\lambda dz_1)u=\v$
and further results}
\smallskip

Let $\hat R$ be the operator defined by formula (2.4)  and let $R_{\lambda}$ be the
operator defined by formula (3.1).
\bigskip

\noindent
{\bf Proposition 4.}
{\sl
Let $\v\in L_{1,1}^{\infty}(V)$ with support in
$V_0=\{z\in V:\ |z_1|\le r_0\}$, where $r_0$ satisfies the condition of
$\S 1$. Then, for
$u=G_{\lambda}\v\buildrel \rm def \over =R_{\lambda}\circ\hat R\v$, where
$\lambda\ne 0$, one has
\smallskip
\item{  i)} \hskip1.3cm$\bar\pa(\pa+\lambda dz_1)u=\v+\bar\lambda d\bar z_1\wedge
{\cal H}_{\lambda}(\hat R\v)$ on $V$;
\smallskip
\item{ ii)}
\vskip-.7cm
$$\eqalign{
&\|u\|_{L^{\infty}(V)}\le {\rm const}(V_0,\tilde p)\cdot
\min\,(1/\sqrt{|\lambda|},1/|\lambda|)\,\|\v\|_{L_{1,1}^{\infty}(V_0)},\ \
\tilde p>2,\cr
&\|\pa u\|_{L^{\tilde p}_{1,0}(V)}\le
{\rm const}(V_0,\tilde p)\,\|\v\|_{L_{1,1}^{\infty}(V_0)},\ \ \tilde p>2.\cr}$$
}
\bigskip

\noindent
{\bf Supplement.}
If we can write $\v=\v_0+\v_1$, where $\v_0\in L^{\infty}_{1,1}(V)$,
${\rm supp}\,\v_0\subset V_0$, and $\v_1=\sum\limits_{l=1}^gic_l\delta(z,a_l)$,
with $a_l\in V_{j(l)}\cap\tilde V_0$, then instead of i)-ii) we have i) and the following conclusion:
\smallskip
{\sl \item{ii)$^{\prime}$}
\vskip-.8cm
$$\eqalign{
&\|u-u(\infty_l)\|_{L^{2+\ep}(V_l)}\le const(V,\ep)\cdot
\min\,(|\lambda|^{-1/2},|\lambda|^{-1})
\bigl(\|\v_0\|_{L^{\infty}_{1,1}(V_0)}+\sum_{j=1}^g|c_j|\bigr),\cr
&\|\pa u\|_{L^{2\pm\ep}_{1,0}(V)}\le const(V,\ep)
\bigl(\|\v_0\|_{L^{\infty}_{1,1}(V_0)}+\sum_{l=1}^g|c_l|\bigr),\cr}$$
where $0<\ep<1/2$.}
\bigskip

\noindent
{\it Proof:}
By Proposition 2 we have
$$f=Fdz_1=\hat R\v\in {\tilde W}^{1,\tilde p}_{1,0}(V)\ \
\forall\tilde p\in (2,\infty),\ \
F\big|_{V_0}\in W^{1,p}(V_0)\ \forall p\in (1,2).$$
Propositions 2 and 3 imply that
$u=R_{\lambda}\circ\hat R\v\in {\tilde W}^{1,\tilde p}(V)$.
Let us now verify statement i) of Proposition 4. From Proposition 3 i) we
obtain
$$\eqalign{
&(\pa+\lambda dz_1)u=(\pa+\lambda dz_1)R_{\lambda}\circ\hat R\v=
\hat R\v+{\cal H}_{\lambda}(\hat R\v),\ \ {\rm where}\cr
&{\cal H}_{\lambda}(\hat R\v)=e_{-\lambda}{\cal H}(e_{\lambda}\hat R\v)
.\cr}\eqno(4.1)$$
From (4.1) and Proposition 2 we obtain
$$\bar\pa(\pa+\lambda dz_1)u=\v+\bar\pa({\cal H}_{\lambda}(\hat R\v))=\v+
\bar\lambda d\bar z_1\wedge {\cal H}_{\lambda}(\hat R\v),$$
where we have used that ${\cal H}(\hat R\v)\in H_{1,0}(\tilde V)$.

\noindent Property 4 ii)  follows from Proposition 3 ii), iii).
The supplement to Proposition 4 follows from the supplement to Proposition 3.
\bigskip

\noindent
{\it Definition}

\noindent We define the Faddeev type Green function for $\bar\pa(\pa+\lambda dz_1)$ on
$V$ as the kernel $g_{\lambda}(z,\xi)$ of the integral operator
$R_{\lambda}\circ\hat R$.
\bigskip

\noindent
{\it Definition}

\noindent Let  $q\in C_{1,1}(\tilde V)$ be a form with  ${\rm supp}\,q$ contained in $V_0$, and
let $g$ denote the genus of $\tilde V$.
The function $\psi(z,\lambda)$, $z\in V$, $\lambda\in\C$, will be called
the Faddeev type function associated with the potential $q$ and the points
$a_1,\ldots,a_g\in  V\b \bar V_0$, if \ $\forall\lambda\in\C\b E$,
where $E$ is compact in $\C$, the
function $\mu=\psi(z,\lambda)e^{-\lambda z_1}$ satisfies the properties:
\vskip-.2cm
$$\eqalign{
&\bar\pa(\pa+\lambda dz_1)\mu={i\over 2}q\mu+i
\sum_{l=1}^gc_l\delta(z,a_l)\ \ {\rm and} \ \
\lim\limits_{\scriptstyle z\to\infty \atop\scriptstyle z\in V_1}
\mu(z,\lambda)=1,\cr
&(\mu-\mu(\infty_j))\big|_{V_j}\in L^{\tilde p}(V_j),\ \ \tilde p>2,
\ j=1,\ldots,d,\cr}$$
where $\delta(z,a_l)$- Dirac measure concentrated in point $a_l$.
\bigskip

Based on the Faddeev type Green function $g_{\lambda}(z,\xi)$, and on
Proposition 4,   we have in [HM] extended the
Novikov reconstruction scheme  from the case $X\subset\C$ to the case of
a bordered Riemann surface $X\subset V$.
\bigskip

\noindent
{\bf Definition}

\noindent Let $\{\omega_j\}$ be an orthonormal basis for the holomorphic 
forms on $\tilde V$.
An effective divisor $\{a_1,\ldots,a_g\}$ on $V$ will be called generic, if
$$\det\bigl[{\omega_j\over d z_1}(a_k)\big|_{j,k=1,2,\ldots,g}\bigr]\ne 0.$$
\bigskip

\noindent
{\bf Lemma.}
{\sl
Let $\{a_j\}$ be a generic divisor on $V$. Put
$$\Delta(\lambda)=\det\bigl[\int\limits_{\xi\in V}\hat R(\delta(\xi,a_j))
\wedge\bar\omega_l(\xi)
e^{\lambda\xi_1-\bar\lambda\bar\xi_1}\bigr]\big|_{j,l=1,2,\ldots,g},$$
where $\hat R$ is the operator from Proposition 2. Then, under the condition 
that
$|a_j|\ge A$, $j=1,2,\ldots,g$, with $A$ large enough,
$\overline{\lim}_{\lambda\to\infty}|\lambda^g\cdot\Delta(\lambda)|<
\infty$,
$\underline{\lim}_{\lambda\to\infty}|\lambda^g\cdot\Delta(\lambda)|>0$
and the set
$$E=\{\lambda\in\C:\ \Delta(\lambda)=0\}\ \
{\rm is\  a\ compact\ nowhere\ dense\ subset\ of}\ \ \C.\eqno(*)$$ }

The following is a corrected version of the main results from [HM]:
\smallskip
\item{1.} Let $X$ be a domain with smooth boundary on $V$ such that
$X\supset\bar V_0$, $\bar X\subset Y\subset V$.
Let $\sigma\in C^{(2)}(V)$, $\sigma>0$ on $V$ and $\sigma=1$ on
$V\b X$.
Let $a_1,\ldots,a_g$ be a generic divisor on $Y\b\bar X$, satisfying condition
$(*)$.
Then for all $\lambda\in\C\b E$ there exists a unique Faddeev type
function $\psi(z,\lambda)=\mu(z,\lambda)e^{\lambda z_1}$ associated with the
potential $q={d\,d^c\sqrt{\sigma}\over \sqrt{\sigma}}$ and the divisor
$\{a_j\}$. Such a function can be found (together with constants $\{c_l\}$)
from the integral equation:
$$\mu(z,\lambda)=1+{i\over 2}\int\limits_{\xi\in X}
g_{\lambda}(z,\xi)\mu(\xi,\lambda)q(\xi)
+i\sum_{l=1}^gc_l(\lambda)g_{\lambda}(z,a_l),\eqno(4.2)$$
where
$$\eqalign{
&{1\over 2}{\cal H}_{\lambda}(\hat R(q \mu))=
\sum_{l=1}^gc_l{\cal H}_{\lambda}(\hat R(z,a_l)),\cr
&\mu(z,\lambda)\to 1,\ \ z\in V_1,\ \ z\to\infty,\cr}\eqno(4.3)$$
$\lambda\in\C\b E$.
\smallskip

The relation (4.3) is equivalent to the system of equations
$$\eqalign{
&2\sum_{l=1}^gc_l(\lambda)e^{\lambda a_{j,1}-\bar\lambda\bar a_{j,1}}
{\bar\omega_k\over d\bar z_1}(a_j)=\cr
&-\int\limits_{z\in X}
e^{\lambda z_1-\bar\lambda\bar z_1}
\bigl({d\,d^c\sqrt{\sigma}\over \sqrt{\sigma}}-
2i\pa\ln\sqrt{\sigma}\wedge\bar\pa\ln\sqrt{\sigma}\bigr)
\mu(z,\lambda){\bar\omega_k\over d\bar z_1}(z),\cr}$$
\item{} where $k=1,\ldots,g$ and
$\{\omega_j\}$ is an orthonormal basis of holomorphic forms on $\tilde V$.
\smallskip

\item{2.} For all $\lambda\in\C\b E$ the restriction of
$\mu=e^{-\lambda z_1}\psi(z,\lambda_1)$ to $bX$ can be found through
Dirichlet-to-Neumann data for $\mu$ on $bX$ by the Fredholm integral equation
$$\mu(z,\lambda)\big|_{bX}+\int\limits_{\xi\in bX}g_{\lambda}(z,\xi)
(\bar\pa\mu(\xi,\lambda)-\bar\pa\mu_0(\xi,\lambda))=
1+i\sum_{j=1}^gc_jg_{\lambda}(z,a_j),\eqno(4.4)$$
where
$$-i\sum_{j=1}^g(a_{j,1})^{-k}c_j=\int\limits_{z\in bX}
z_1^{-k}(\pa+\lambda dz_1)\mu=
0,\quad k=2,\ldots,g+1,\eqno(4.5)$$
and $\mu_0$ is the solution of the Dirichlet problem
$$\bar\pa(\pa+\lambda dz_1)\mu_0\big|_X=0,\ \  \mu_0\big|_{bX}=\mu\big|_{bX}.
$$
The parameters $\{a_{j,1}\}$ (the first coordinates of $\{a_j\}$) are supposed to be
mutually different.
\smallskip

The equations (4.4), (4.5) are solvable simultaneously with (4.2), (4.3).
\smallskip

The relations (4.5) are equivalent to the equality
$$\bar\pa(\pa+\lambda dz_1)\mu\big|_{V\b X}=i\sum_{j=1}^gc_j\delta(z,a_j).$$
\item{3.} The Faddeev type function
$\mu=\psi(z,\lambda)e^{-\lambda z_1}$ satisfies the
Bers--Vekua type $\bar\pa$-equation with respect to $\lambda\in\C\b E$
$${\pa\mu(z,\lambda)\over \pa\bar\lambda}=b(\lambda)\bar\mu(z,\lambda)
e^{\bar\lambda\bar z_1-\lambda z_1},\eqno(4.6)$$
where
$$b(\lambda)\buildrel \rm def \over =
\lim\limits_{\scriptstyle z\to\infty \atop\scriptstyle z\in V_l}
{\bar z_1\over \bar\lambda}e^{\lambda z_1-\bar\lambda \bar z_1}
{\pa\mu\over \pa\bar z_1}(z,\lambda)\bigg/
\lim\limits_{\scriptstyle z\to\infty \atop\scriptstyle z\in V_l}
\overline{\mu(z,\lambda)},$$
with $l=1,\ldots,d$.
The function $b(\lambda)$, referred to as nonphysical scattering data, can be found
by (4.6) through $\mu\big|_{bX}$.
\smallskip

In addition, the following important formulas for the data $b(\lambda)$ are
valid
$$d\cdot\bar\lambda\cdot b(\lambda)=-{1\over 2\pi i}\int\limits_{z\in bY}
e^{\lambda z_1-\bar\lambda\bar z_1}\bar\pa\mu=
{1\over 2\pi i}\int\limits_{z\in X}{i\over 2}
e^{\lambda z_1-\bar\lambda\bar z_1}q\mu+i\sum_{j=1}^gc_j
e^{\lambda a_{j,1}-\bar\lambda\bar a_{j,1}},\eqno(4.7)$$
\item{} where $\lambda\in\C\b E$.
\smallskip

\item{} On the basis of (4.3), (4.7) and Proposition 3, one can derive the estimate
$$|\lambda\cdot b(\lambda)|\le const (V,\sigma)
(1+|\lambda|)^{-g}
|\Delta(\lambda)|^{-1},\ \ \lambda\in\C\b E.
\eqno(4.8)$$

\item{4.}
Let us suppose now that the divisor $\{a_1,\ldots,a_g\}$ on $Y\b X$
is such that the exceptional compact $E$ in $\C$ consists of isolated points
$\lambda_1,\ldots,\lambda_N$ and
$$|\Delta(\lambda)|\ge {\rm const}(V) {\rm dist}(\lambda,E)\ \ {\rm if}\ \
{\rm dist}(\lambda,E)\le {\rm const}.\eqno(4.9)$$
Then the reconstruction procedure for $\mu\big|_{X\times\C}$ and
$\sigma\big|_X$ through scattering data $b\big|_{\C}$ can be done in the
following way.
\smallskip

The relations (4.2), (4.3), combined with the inequalities (4.8), (4.9), imply that the
\item{}$\bar\pa$-equation (4.6) can be replaced by  the singular integral equation:
$$\eqalign{
&(\mu-1)+{1\over 2\pi i}\lim\limits_{\delta\to 0}
\int\limits_{\C\b\cup\{|\xi-\lambda_l|\le\delta\}}
b(\xi)
e^{\bar\xi\bar z_1-\xi z_1}\overline{(\mu-1)}
{d\bar\xi\wedge d\xi\over {\xi-\lambda}}+
{1\over 2\pi i}\sum_{l=1}^N{\mu_l\over {\lambda_l-\lambda}}=\cr
&-{1\over 2\pi i}\int\limits_{\C}b(\xi)
e^{\bar\xi\bar z_1-\xi z_1}
{d\bar\xi\wedge d\xi\over {\xi-\lambda}},\ \ {\rm where}\cr
&\mu_l=\lim\limits_{\delta\to 0}\int\limits_{|\xi-\lambda_l|\le\delta}
b\bar\mu e^{\bar\xi\bar z_1-\xi z_1}d\bar\xi\wedge d\xi=
\lim\limits_{\delta\to 0}\int\limits_{|\xi-\lambda_l|=\delta}
\mu d\xi=O_z(1),\cr
&l=1,2,\ldots,N,\ \ \lambda\in\C\b E.\cr}\eqno(4.10)$$
This equation is of Fredholm--Noether type in the  space of functions
$$\lambda\mapsto (\mu(\cdot,\lambda)-1):\ |\mu-1|\cdot |\Delta(\lambda)|
(1+|\lambda|)\in L^{\tilde p}(\C),\ \ \tilde p>2.$$
In contrast to the planar case,
when $d=1$, $g=0$, equation (4.10) does not necessarily have a unique solution.
This makes it possible for almost all $z_1\in\C$ to find a basis of
independent solutions of (4.10)
$$\lambda\mapsto\mu_k(z_1,\lambda),\ \ k=1,2,\ldots,\tilde d,\ \
\lambda\in\C,\ \ \tilde d\ge d.$$
Put
$$\mu(z_1,z_2,\lambda)=\mu(z_1,z_{2,j}(z_1),\lambda)=\sum_{k=1}^{\tilde d}
\gamma_{j,k}(z_1)\mu_k(z_1,\lambda),$$
where $(z_1,z_2)=(z_1,z_{2,j}(z_1))\in V$, $j=1,2,\ldots,\tilde d$.
The condition for the form $\mu^{-1}\bar\pa(\pa+\lambda dz_1)\mu$ to be
independent of $\lambda$ allows us to find (maybe not uniquely) the coefficients
$\gamma_{j,k}(z)$ in the
expression for $\mu(z_1,z_2,\lambda)$.   The equalities
$${i\over 2}{d\,d^c\sqrt{\sigma}\over \sqrt{\sigma}}\big|_X=
q\big|_X=\mu^{-1}\bar\pa(\pa+\lambda dz_1)\mu\big|_X$$
finally permit us to find all $q$ and $\sigma$ with given scattering
data $b\big|_{\C}$.
\bigskip

\noindent The uniqueness of the reconstruction of $\mu\big|_{X\times\C}$ and
$\sigma\big|_X$
from the data $b$ on $\C\setminus E$ is plausible but still unknown.
Nevertheless, the uniqueness of the reconstruction of $\sigma\big|_X$ from
Dirichlet-to-Neumann data of the equation $d(\sigma d^cU)\big|_X=0$ can be
proved by the above procedure using Dirichlet-to-Neumann data not just for a single function,
but for a family of Faddeev type functions depending on a parameter $\theta$:
$$\eqalign{
&\psi_{\theta}(z,\lambda)=e^{\lambda(z_1+\theta z_2)}
\mu_{\theta}(z_1,z_2,\lambda),\ \ {\rm where}\cr
&\bar\pa(\pa+\lambda(dz_1+\theta dz_2))\mu_{\theta}={i\over 2}q\mu_{\theta}+i
\sum_{l=1}^gc_l\delta(z,a_l)\ \ {\rm and} \ \
\lim\limits_{\scriptstyle z\to\infty \atop\scriptstyle z\in V_1}
\mu_{\theta}(z,\lambda)=1,\cr
&(\mu_{\theta}-\mu_{\theta}(\infty_j))\big|_{V_j}
\in L^{\tilde p}(V_j),\ \ \tilde p>2,\ \
\lambda\in\C\b E_{\theta},\ \ j=1,\ldots,d.\cr}$$

For the reconstruction of $\sigma\big|_X$ it is in fact sufficient to use data
$\psi_{\theta}(z,\lambda)\big|_{bX\times\C}$ for at most $d$ different values of
the parameter $\theta$.

\vskip 4 mm
{\bf References}

\item{[ BC1]} Beals R., Coifman R., Multidimensional inverse scattering and
nonlinear partial differential equations, Proc.Symp. Pure Math. {\bf 43}
(1985), A.M.S. Providence, Rhode Island, 45-70
\item{[ BC2]} Beals R., Coifman R., The spectral problem for the
Davey-Stewartson  and Ishimori hierarchies, In: "Nonlinear Evolution Equations
: Integrability and Spectral Methodes", Proc. Workshop, Como, Italy 1988,
Proc. Nonlinear Sci., 15-23, 1990
\item{[  F1]} Faddeev L.D., Increasing solutions of the Schr\"odinger
equation, Dokl.Akad.Nauk SSSR {\bf 165} (1965), 514-517 (in Russian);
Sov.Phys.Dokl. {\bf 10} (1966), 1033-1035
\item{[  F2]} Faddeev L.D., The inverse problem in the quantum theory of
scattering, II, Current Problems in Math., {\bf 3}, 93-180, VINITI, Moscow,
1974 (in Russian); J.Soviet Math. {\bf 5} (1976), 334-396
\item{[  GN]} Grinevich P.G., Novikov S.P, Two-dimensional "inverse
scattering problem" for negative energies and generalized analytic
functions, Funct.Anal. and Appl., {\bf 22} (1988), 19-27
\item{[   H]} Henkin G.M., Uniform estimate for a solution of the $\bar\pa$-
equation in Weil domain, Uspekhi Mat.Nauk {\bf 26} (1971), 211-212
\item{[  HM]} Henkin G.M., Michel V., Inverse conductivity  problem
on Riemann surfaces,
\item{      }J.Geom.Anal. (2008) {\bf 18}(4), 1033-1052
\item{[  HP]} Henkin G.M., Polyakov P.L., Homotopy formulas for the $\bar\pa$-
operator on ${\C}P^n$ and the Radon-Penrose transform, Math. USSR Izvestiya
{\bf 28} (1987), 555-587
\item{[  Ho]} Hodge W., The theory and applications of harmonic integrals,
Cambridge Univ. Press, 1952
\item{[ N1]} Novikov, R., Multidimensional inverse spectral problem for the
equation $-\Delta\psi+(v-Eu)\psi=0$, Funkt.Anal. i Pril. {\bf 22} (1988),
11-22 (in Russian); Funct.Anal and Appl. {\bf 22} (1988), 263-278
\item{[ N2]} Novikov, R., The inverse scattering problem on a fixed energy
level for the two-dimensional Schr\"odinger operator, J.Funct.Anal. {\bf 103}
(1992), 409-463
\item{[  Na]} Nachman A., Global uniqueness for a two-dimensional inverse
boundary problem, Ann. of Math. {\bf 143} (1996), 71-96
\item{[  P1]} Pompeiu D., Sur la repr\'esentation des fonctions analytiques
par
des int\'egrales d\'efinies, C.R.Acad.Sc.Paris {\bf 149} (1909), 1355-1357
\item{[P2]} Pompeiu D., Sur une classe de fonctions d'une variable complexe
et sur certaines \'equations int\'egrales, Rend. del. Circolo Matem. di
Palermo {\bf 35} (1913),277-281
\item{[ Po]} Polyakov P., The Cauchy-Weil formula for differential forms,
Mat.Sb. {\bf 85} (1971), 388-402
\item{[   S]} Serre J.-P., Un th\'eor\`eme de dualit\'e, Comm.Math.Helv.
{\bf 29} (1955), 9-26
\item{[   V]} Vekua I.N., Generalized analytic functions, Pergamon Press, 1962
\item{[   W]} Weil A., Vari\'et\'es k\"ahleriennes, Hermann, Paris, 1957

\vskip 4 mm
\line{\hfill\box3}

\end